\documentclass[12pt]{amsart}
\font\bbbld=msbm10 scaled\magstephalf

\newcommand{\bM}{\bar{M}}

\newcommand{\bfR}{\hbox{\bbbld R}}

\newcommand{\cC}{\mathcal{C}}

\newcommand{\tr}{\mbox{tr}}

\newcommand{\tu}{\tilde{u}}
\newcommand{\tvarphi}{\tilde{\varphi}}

\newcommand{\ol}{\overline}
\newcommand{\ul}{\underline}

\newtheorem{theorem}{Theorem}[section]

 \theoremstyle{definition}

\theoremstyle{remark}
\newtheorem{remark}[theorem]{Remark}

\numberwithin{equation}{section}

\begin{document}

\setlength{\baselineskip}{1.2\baselineskip}

\title[The Neumann problem fo Fully Nonlinear Elliptic Equations]
{On estimates for fully nonlinear elliptic equations with Neumann boundary conditions on Riemannian Manifolds}
\author{Bo Guan}
\address{Department of Mathematics, Ohio State University,
         Columbus, OH 43210}
\email{guan@math.osu.edu}
\author{Ni Xiang}
\address{Department of Mathematics and Statistics, HuBei University,
WuHan, P.R.China}
\email{nixiang@hubu.edu.cn}
\thanks{Research supported in part by NSF grants (Guan), CSC 201708420065 and Hubei Key
Laboratory of Applied Mathematics in Hubei University
(Xiang).}

\begin{abstract}
We derive gradient and second order {\em a priori} estimates for solutions of 
the Neumann problem for a general class of fully nonlinear elliptic equations 
on compact Riemannian manifolds with boundary.  
These estimates yield regularity and existence results.


\end{abstract}

\maketitle
\bigskip

\section{Introduction}

\medskip

 Let $(\bM^n, g)$ be a compact Riemannian manifold
of dimension $n \geq 2$ with smooth boundary $\partial M$,
and let $M$, $\nabla$ denote the interior of $\bM$ and the Levi-Civita connection of $g$, respectively.
For a function $u \in C^2 (\bM)$ and a 2-tensor $A$
on $M$, let $\nabla^2 u$ denote the Hessian of $u$
and $\lambda [A] = (\lambda_1, \cdots, \lambda_n)$
the eigenvalues of $A$ with respect to the metric $g$.
In this paper we are concerned with fully nonlinear elliptic equations of the form
\begin{equation}
 \label{021901}
 f ( \lambda [\nabla^2 u + \chi]) = \psi
     \;\; \mbox{in $\bM$}
\end{equation}
with Neumann boundary data
\begin{equation}
\label{021906}
\nabla_\nu u = \varphi (x, u)
\;\; \mbox{on $\partial M$},
\end{equation}
where $\nu$ denotes the interior unit normal to $\partial M$,
and $f$ is a smooth symmetric
function of $n$ variables
defined in a symmetric open and convex cone $\Gamma$ in $\bfR^n$ containing $\Gamma_n$
with vertex at the origin, where
\begin{equation}
\label{021902}
\Gamma_n \equiv \{\lambda \in \bfR^n:
\mbox{all $\lambda_i>0$}\}.
\end{equation}

Following \cite{CNS85} we assume $f$ to satisfy the
fundamental structure conditions
\begin{equation}
\label{021903}
f_i 
\equiv \frac {\partial f} {\partial \lambda_i}>0 \; \mbox{in} \; \Gamma, \; 1 \leq i \leq n,
\end{equation}
\begin{equation}
\label{021904}
\mbox {$f$ is a concave function in $\Gamma$},
\end{equation}
and
\begin{equation}
\label{021905}
\sup_{\partial \Gamma} f < \inf _{M} \psi
\end{equation}
where
\[  \sup_{\partial \Gamma} f\equiv \sup_{\lambda_0 \in \partial \Gamma} \lim \sup_{\lambda \rightarrow \lambda_0} f (\lambda). \]


A function $u \in C^2(M)$ is called {\em admissible} if $ \lambda [\nabla^2 u + \chi] \in \Gamma$. Condition (\ref{021903}) ensures that (\ref{021901}) is elliptic for admissible solution $u \in C^2(M)$, while (\ref{021904}) implies the function $F$ defined by $F(A) = f (\lambda[A])$ is concave for
$A \in \mathcal{S}^{n \times n}$
with $\lambda[A] \in \Gamma$, where $\mathcal{S} ^ {n \times n}$ denotes the set of $2$-tensors on $M$. Condition (\ref{021905}) prevents equation (\ref{021901}) from being degenerate; see \cite{CNS85}.

The most typical examples of form (\ref{021901}) are given by $f = \sigma^{1/k}_k$ and $f = (\sigma_k / \sigma_l) ^ {1/(k-l)}$, $1 \leq l< k \leq n$, defined on the cone
\begin{equation}
\Gamma_k = \{\lambda \in \bfR^n: \sigma_j (\lambda)>0 \  \mbox{for}\  1 \leq j \leq k\},
\end{equation}
where $\sigma_k$ is the $k$-th elementary symmetric function
\begin{equation}
\sigma_k (\lambda) = \sum_{i_1 < \cdots < i_k} \lambda_{i_1} \cdots \lambda_{i_k}.
\end{equation}
There are other interesting functions satisfying
\eqref{021903} and \eqref{021904} which
naturally arise from important geometric problems.

The Neumann problem in uniformly convex domains in $\bfR^n$ was treated 
by Lions-Trudinger-Urbas  \cite{LTU}  for Monge-Amp\`ere equation 
and recently by Ma-Qiu \cite{MQ} for Hessian equations corresponding to 
$f = \sigma^{1/k}_k$.
Urbas \cite{Ur95, Ur96} studied the oblique boundary value problems 
for Hessian and curvature equations in two dimensions. 
Meanwhile, the Dirichlet problem has received much more attention and 
extensive study since  work of Ivochkina \cite{Iv80} and 
Caffarelli-Nirenberg-Spruck~\cite{CNS85}; 
see e.g.  \cite{CW01, Guan14, Guan, Tr95, Wang94}, etc.

Our primary goal in this paper is to establish {\em a priori} estimates for admissible solutions of
problem~\eqref{021901}-\eqref{021906}.
For this we recall some notions and results from \cite{Guan14}.

For $\sigma > \sup_{ \partial \Gamma} f$, set
$\Gamma ^ {\sigma} = \{\lambda \in \Gamma: f (\lambda) > \sigma\}$.
By conditions~(\ref{021903}) and (\ref{021904}),
the boundary of $\Gamma^{\sigma}$
\[ \partial \Gamma ^ {\sigma} = \{\lambda \in \Gamma: f (\lambda) = \sigma\}, \]
which is a level hypersurface of $f$,
is smooth and convex.
Define for $\mu \in \Gamma$
 \[ S^{\sigma}_{\mu} = \{\lambda \in \partial
\Gamma^{\sigma}: \nu_{\lambda}\cdot (\mu-\nu)\leq 0 \}\]
and
\[ \mathcal{C} _ {\sigma}^ + = \{\mu \in \Gamma: S_{\mu} ^{\sigma} \mbox { is compact}\}.
 \]
It was shown in \cite{Guan14} that
$\mathcal{C}_{\sigma}^+$ is open.
We call $\partial \mathcal{C}_{\sigma}^+$ {\em the
tangent cone at infinity} of $\Gamma^{\sigma}$.
The following assumption plays key roles in our results:
there exist an admissible  function $v \in C^2(\bM)$ satisfying
\begin{equation}
\label{032201}
\begin{aligned}
&\lambda [\nabla^2 v + \chi] \in \mathcal{C}_{\psi}^+.
\end{aligned}
\end{equation}

Throughout of this paper $\psi$ and $\varphi$ are assumed to be smooth functions, and let
$u \in C^4 (M) \cap C^3 (\bM)$ be an admissible solution of the Neumann problem \eqref{021901}-\eqref{021906}.
Our first main result concerns the global second derivative estimates and may be stated as follows.

\begin{theorem}
\label{mainth1}
Under conditions~\eqref{021903}-\eqref{021905} and
\eqref{032201}, there exists a constant $C$
depending
$|u|_{C^1 (\bM)}$ and other known data such that
\begin{equation}
\label{hess-a10}
\max_{\bar{M}} |\nabla^2 u| \leq C \Big(1
+\max_{\partial M}|\nabla^2 u|\Big).
\end{equation}
If moreover,
\[ \varphi = a (x) u + b (x), \]
where $a(x)$ and $b(x)$ are smooth functions,
  then
\begin{equation}
\label{hess-a10a}
\max_{\bar{M}} |\nabla^2 u| \leq
    C_1 
  \Big(1 + \max_{\bar M} |\nabla u|^2
            +\max_{\partial M}|\nabla^2 u|\Big),
\end{equation}
where $C_1$ is
independent of $|\nabla u|_{C^0 (\bM)}$.
\end{theorem}

Our next result concerns the gradient estimates.

\begin{theorem}
\label{gx-thm-I20}
Assume in addition to \eqref{021903}-\eqref{021905} and \eqref{032201} that
\begin{equation}
\label{082301}
 \sum f_i (\lambda) \lambda_i
   \geq - \omega_f (|\lambda^-|) \sum f_i
   \; {\mbox { in $\Gamma ( \psi ) $}}
\end{equation}
when $|\lambda|$ is sufficiently large,
where  $\lambda^- = (\lambda_1^-, \ldots, \lambda_n^-)$,
\[ \Gamma(\psi) := \Gamma \cap
           \{ \inf _M \psi \leq f \leq \sup _M \psi \}, \]
and $\omega_f \geq 0$ is a nondecreasing function
satisfying the sublinear growth condition
\begin{equation}
\label{gx-I110}
\lim_{t \rightarrow + \infty} \frac{\omega_f (t)}{t} = 0.
\end{equation}
Then
\begin{equation}
\label{082302}
\max_ {\bM} |\nabla u| \leq C.
\end{equation}
for some constant $C$ depending on $|u|_{C^0 (\bM)}$
and other known data.
\end{theorem}

\begin{remark}
It is not clear to the authors if \eqref{082301} in fact always holds for some $\omega_f$ satisfying
\eqref{gx-I110}; obviously it does for
$\omega_f (t) = t$.
\end{remark}

Turning to the boundary estimates for second order derivatives, we need to strengthen 
our assumptions and impose restrictions to the underlying manifold and its boundary.

\begin{theorem}
\label{gx-thm-I30}
Suppose that $(M^n, g)$ is locally conformally flat
near boundary and $\partial M$ is umbilic.
Assume in addition to \eqref{021903}-\eqref{021905} and \eqref{032201}  
that
\begin{equation}
\label{gx-I110a}
\sum f_i \lambda_i \geq 0 \;\; \mbox{in $\Gamma (\psi)$}
\end{equation}
and that
the function $v$ 
satisfies
\begin{equation}
\label{032201a}
 \nabla_{\nu} v \geq \sup_{a \leq t \leq b} \varphi (x, t) + \epsilon_0 \; \mbox{on $\partial M$}
\end{equation}
for some constant $\epsilon_0 > 0$ where
$a = \inf_M u$, $b = \sup_M u$.
Suppose furthermore that there exists a fucntion
$\ul u \in C^2 (\bM)$ satisfying
\begin{equation}
\label{032301} 
\left\{
\begin{aligned}
& \lambda [\nabla^2 \ul u] \in \cC^+_{\psi}
      \;\; \mbox {in $M$},\\
& f (\lambda [ \nabla^2 \underline{u} ]) \geq \psi
      \;\; \mbox {in $M$},\\
& \underline{u} = 0\;\; \mbox{on $\partial M$}.
\end{aligned}
\right.
\end{equation}
Then
\begin{equation}
\label{hess-a10b}
\max_{\partial M} |\nabla^2 u| \leq C
\end{equation}
for some constant $C$ depending
$|u|_{C^1 (\bM)}$ and other known data.
If moreover,
\[ \varphi = a (x) u + b (x) \]
 then
\begin{equation}
\label{hess-a10b'}
\max_{\partial M} |\nabla^2 u| \leq
 C_1 \Big(1 + \max_{\bar M} |\nabla u|^2\Big)
\end{equation}
where $C_1$ is
independent of $|\nabla u|_{C^0 (\bM)}$.
\end{theorem}

More generally, we say
$(\bM^n, g)$ is $\Gamma$-{\em admissible} if there exists a function $\ul u \in C^2 (\bM)$ with
$\lambda [\nabla^2 \ul u] \in \Gamma$ in $\bM$ and $\ul u = 0$ on $\partial M$.
In particular, $\ul u$ is subharmonic, i.e.
$\Delta \ul u \geq 0$ and therefore $\ul u \leq 0$ in $\bM$.

\begin{remark}
If the function $\varphi$ in the Neumann condition~\eqref{021906}
is independent of $u$, in place of $|u|_{C^0 (\bM)}$
the constants $C$ and $C_1$ in the above estimates will only depnds on
\[ \omega (u) := \sup_M u - \inf_M u. \]
\end{remark}

 \begin{remark}
 For $f = \sigma^{1/k}_k$, $k \geq 2$, $\partial \Gamma^ {\sigma}$ is strictly convex and
$\mathcal{C}_ {\sigma}^ + = \Gamma_k$ for all $\sigma>0$.This follows from the property
\begin{equation}
\label{3I-50}
\lim_{R \rightarrow}
 f (\lambda_1, \ldots, \lambda_{n-1}, \lambda_n + R)
  = + \infty.
\end{equation}
\end{remark}

Based on the estimates in Theorems~\ref{mainth1}, \ref{gx-thm-I20}, \ref{gx-thm-I30} and the
Evans-Krylov theorem due to
Lieberman and Trudinger~\cite{LT}
we obtain the following existence result on the Neumann problem~(\ref{021901})-(\ref{021906}).

\begin{theorem}
\label{3I-thm2}
Let $(M^n, g) $ be a smooth compact Riemmannian manifold with smooth umbilic boundary $\partial M$.
Suppose that $M$ is locally conformally flat near boundary and that for some constant $\gamma_0 > 0$,
\begin{equation}
\label{gx-I140}
\nabla_{\nu} u \geq \gamma_0.
\end{equation}
Assume that conditions~\eqref{021903}-\eqref{021905}, 
\eqref{032201} and 
\eqref{gx-I110a}-\eqref{032301} hold 
with $\ul u$ furthermore satisfying
\begin{equation}
\label{gx-I155}
f (\lambda [\nabla^2\ul u]) \geq \psi
\;\; \mbox{in $\bM$}
\end{equation}
and $a$, $b$ in \eqref{032201a} given by
\begin{equation}
\label{gx-I145}
 a = \inf_{M} \ul u
    + \frac{1}{\gamma_0} \inf_{x \in M} \varphi (x, 0),
    \;\;
b = \sup_{M} h + \frac{1}{\gamma_0} \sup_{x \in M} \varphi (x, 0),
\end{equation}
respectively, where $h \in C^{\infty} (\bM)$ satisfies
$\Delta h + \tr \, \chi = 0$ in $\bM$ and $h = 0$ on $\partial M$.
Then the Neumann problem~\eqref{021901}-\eqref{021906} admits a unique admissible solution
$u \in C^{\infty} (\bM)$.
\end{theorem}

Indeed, by \eqref{021906} and \eqref{gx-I140}
one immediately derives the lower and upper bounds by the maximum principle,
\begin{equation}
\label{gx-I165}
 \inf_{M} \ul u
    + \frac{1}{\gamma_0} \inf_{x \in M} \varphi (x, 0)
    \leq u \leq
\sup_{M} h + \frac{1}{\gamma_0} \sup_{x \in M} \varphi (x, 0).
\end{equation}
Consequently, one can apply the continuity method and the classical Schauder theory to
prove the existence of an smooth admissible solution.

Equations of form (\ref{021901}), especially when $\chi$ and $\psi$ are allowed 
to depend on $u$ and $\nabla u$, naturally appear in  interesting geometric problems 
such as Minkowski problem and the Christoffel-Minkowski problem. 
The Neumann problem for  fully nonlinear equations on manifolds arises in the study 
of fully nonlinear versions of Yamabe problem on manifolds with boundary; we
refer the reader to papers  \cite{Vi00}, \cite{LiLi2006}  \cite{JLL2007}  \cite{Ch07} and references therein.

The rest of this paper is organized as follows. In Section~\ref{gradient}, we shall prove gradient bounds for the admissible solutions.
In Section~\ref{global} and Section ~\ref{B}, we shall derive the global and boundary estimates for second order derivative respectively.

\bigskip

\section{Gradient estimates}\label{gradient}
\label{3I-G}
\setcounter{equation}{0}

\medskip

In this section we derive the gradient estimates.
In order to construct test functions,
we shall first extend data on boundary to the whole manifold, which will also be used in later sections.

Let $d$ be the distance function to $\partial M$.
Since $\partial M$ is smooth, $d$
is smooth in a neighborhood
$M_{\delta_0} := \{x \in M: d (x) < \delta_0\}$
of $\partial M$ for $\delta_0 > 0$ sufficiently small.
 We modify $d$ outside $M_{\delta_0/2}$ into a
smooth function such that  
\[ (\varphi_u + |\varphi_{uu}| + |\varphi_{uuu}|) d \leq \epsilon_0 \]
 for some sufficiently small constant $\epsilon \in (0, \delta_0/2]$
 so in particular $2 \varphi_u d \leq 1$; 
 this will be used in the gradient estimates in this section. 
It is possible to do so by making $\delta_0$ smaller if necessary.
So $\delta_0$ may depend on upper bounds
of $\varphi_u + |\varphi_{uu}| + |\varphi_{uuu}|$.

Note that $\nu = \nabla d$ on $\partial \Omega$
$\nabla d$ denote the gradient of $d$. We may assume
that $\nu$ has been extended to $\bM$ by
$\nu = \nabla d$.
We shall also assume the function $\varphi$ has been smoothly extended to $\bM \times \bfR$,
and write $\tvarphi (x) = \varphi (x, u (x))$.

Let $\tu = u - \tvarphi d$,
$w = 1 + |\nabla \tu|^2$ and
$\phi  = A - \tu - d$ where $A$ is a constant chosen suffieciently large so that $1 \leq \phi \leq 2 A$ in $\bM$.
We assume the function $w \phi^{-1}$ attains a maximum at a point $x_0 \in \bM$.

We first consider the case $x_0 \in \partial M$.
Note that $\nabla_{\nu} \tu = 0$ on $\partial M$ by the boundary condition \eqref{021906}. We may choose orthonormal local frames $e_1, \ldots, e_n$ about $x_0$ so that
$e_n = \nu$ and
$\nabla_1 \tu (x_0) = |\nabla \tu (x_0)| $.
Note that
$d = 0$, $\nabla_1 d = 0$,
$g (\nabla_{n} e_1, e_1) = 0$,
and $\nabla_k \tu = 0$ for $k > 1$
at $x_0$.  We have
$\nabla_1 \tu = \nabla_1 u$,
\[ \nabla_n (\nabla_1 \tu)
  = \nabla_{n1} \tu = \nabla_{1n} \tu
  = \nabla_1 (\nabla_n \tu)
  = - \tvarphi \nabla_{1n} d \]
and consequently,
\begin{equation}
\label{3I-G05b}
\begin{aligned}
  0 \geq \,& \phi \nabla_n w
                      + w \nabla_n (\tu +d) \\
   = \,& \phi \nabla_1 \tu \nabla_n (\nabla_1 \tu)
       + w \\ 
  = \,& - \phi \tilde{\varphi} \nabla_1 u \nabla_{1n} d + w. \\
\end{aligned}
\end{equation}
We derive a bound $w (x_0) \leq C$.

Suppose now that $x_0 \in M$.
We choose
orthonormal local frames about $x_0$ such that
$\nabla_j e_j = 0$, and $U_{ij} := \nabla_{ij} u + \chi_{ij}$ is diagonal at $x_0$.
\begin{equation}
\label{3I-G05a}
   \frac{\nabla_j w}{w} + \frac{\nabla_j \tu+ \nabla_j d}{\phi} = 0
\end{equation}
and
 \begin{equation}
\label{3I-G15a}
\phi F^{ii} \nabla_{ii} w
  + w F^{ii} (\nabla_{ii} \tu + \nabla_{ii} d) \leq 0.
\end{equation}
We calculate
$\nabla_j w = 2 \nabla_k \tu \nabla_{jk} \tu$,
\[ \nabla_{ij} w = 2 \nabla_{ik} \tu \nabla_{jk} \tu
 + 2 \nabla_{k} \tu \nabla_{ijk} \tu, \]
\[ \nabla_k \tu = (1 - \varphi_u d) \nabla_k u -
d \nabla_k \varphi - \tvarphi \nabla_k d, \]
\begin{equation}
\label{3I-G155a}
\begin{aligned}
\nabla_{jk} \tu = \,& (1 - \varphi_u d) \nabla_{jk} u
          -  (\nabla_j u \nabla_k d +
     \nabla_k u \nabla_j d) \varphi_u - d \varphi_{uu} \nabla_j u \nabla_ku \\
     \,& -(\nabla_k \varphi_u \nabla_j u + \nabla_j \varphi_u \nabla_k u) d -
     (\nabla_j d \nabla_k \varphi + \nabla_j \varphi \nabla_k d) \\
      \,&- (d \nabla_{jk} \varphi
    + \tilde{\varphi} \nabla_{jk} d) - (1 - \varphi_u d) \chi_{jk} \\
\end{aligned}
\end{equation}
and
\begin{equation}
\label{3I-G155b}
\begin{aligned}
\nabla_{ijk} \tu = \,&
   (1 - \varphi_u d) \nabla_{ijk} u
          -  (\nabla_{ij} u \nabla_k d +
     \nabla_{ik} u \nabla_j d) \varphi_u - d \varphi_{uuu} \nabla_i u \nabla_j u \nabla_k u
     \\
     \,& - d \varphi_{uu} (\nabla_i u \nabla_{jk} u + \nabla_k u \nabla_{ij} u + \nabla_j u \nabla_{ik} u) \\
     \,& - ( d \nabla_i \varphi_u + \varphi _u \nabla_i d) \nabla_{jk} u
     - \varphi_u (\nabla_j d \nabla_{ik} u + \nabla _k d \nabla_{ij}u )\\
     \,& - d(\nabla_k \varphi_u \nabla_{ij}u + \nabla_j \varphi_u \nabla_{ik} u) + O (w). 
\end{aligned}
\end{equation}

Write equation~\eqref{021901} in the form
\begin{equation}
\label{3I-10'}
F (U) := f (\lambda [U]) = \psi
\end{equation}
for $U = \nabla^2 u + \chi$, and denote\begin{equation}
F^{ij} = \frac {\partial F} {\partial U_{ij}} (U), \;
 F^{ij, kl} = \frac {\partial^2 F} {\partial U_{ij} \partial U_{kl}} (U).
\end{equation}
The matrix $\{F^{ij}\}$ is positive definite by assumption (\ref{021903}) with eigenvalues
$f_1, \cdots, f_n$, and
$F$ is a concave function by (\ref{021904}); see \cite{CNS85}.
Moreover, the following identities hold
for $\lambda [U] = (\lambda_1, \cdots, \lambda_n)$,
\begin{equation}
F^{ij} U_{ij} = \sum f_i \lambda_i,
\end{equation}
and
\begin{equation}
F^{ij} U_{ik} U_{kj} = \sum f_i \lambda_i^2.
\end{equation}

Differentiate equation~\eqref{3I-10'}, we obtain
\begin{equation}
\label{hess-a60}
F^{ij} \nabla_k U_{ij} =  \nabla_k \psi,
\;\;\; \mbox{for all $k$}.
\end{equation}
Note that $\{F^{ij}\}$ is diagonal at $x_0$ since so is $U_{ij}$.
It follows that
\begin{equation}
\label{3I-G25b}
\begin{aligned}
 F^{ii} \nabla_{ii} \tu
 \geq \,& (1 - \varphi_u d) F^{ii} U_{ii}
                - C d w \sum F^{ii}
\end{aligned}
\end{equation}
and, by Schwarz inequality,
\begin{equation}
\label{3I-G25a}
\begin{aligned}
 F^{ii} \nabla_{ii} w
 \geq \,& 2(1 - \varphi_u d)^2 F^{ii} U_{ii}^2 - C |\varphi_{uuu}| d w^2 \sum F^{ii}
 - C d w F^{ii} |U_{ii}|\\
 \,& - C \varphi_{uu}^2 d^2 w^2 \sum F^{ii} - C w^{\frac{3}{2}} \sum F^{ii} - C \sqrt{w} F^{ii}|U_{ii}|\\
               \,& - C w \sum F^{ii}
- 2 (1 - \varphi_u d)^2 |\nabla \psi| \sqrt{w}.
\end{aligned}
\end{equation}

Assume $|\nabla \tu| \geq |\nabla u|/4$ and
let $I = \{i: n |\nabla_i \tu| \geq |\nabla \tu|\}$.
We see that $I \neq \emptyset$
and by \eqref{3I-G05a} for $i \in I$,
\[ (1 - \varphi_u d) U_{ii}
   \leq - \frac{w}{2 \phi} + \frac{w}{|\nabla_i \tu|}
+ C  \leq - \frac{w}{8 A} + C  \equiv - K .\]
We shall assume $w$ to be sufficiently large,
and in particular $K \geq \frac{w}{16 A}$
(otherwise we are done)
so that
we may apply Theorem~2.17 (or Theorem~2.18
more directly) in \cite{Guan14} to
${\mu} = \lambda (\nabla^2 v (x_0)+  \chi (x_0))$ and
$\lambda = \lambda (\nabla^2 u (x_0) + \chi (x_0))$
to derive
\[ - \sum_{U_{ii} < 0} F^{ii} U_{ii}
  \geq F^{ii} (\nabla_{ii} v - \nabla_{ii} u)
           - F^{ii} (\nabla_{ii} v + \chi_{ii})
      \geq \varepsilon - C \sum F^{ii}. \]

Let $J = \{i: U_{ii} \leq - K\}$ so $I \subseteq J$. Clearly,
\begin{equation}
\label{3I-G35a}
 F^{ii} U_{ii}^2 \geq K^2 \sum_J F^{ii}
                    \geq \frac{K^2}{n} \sum F^{ii}.
\end{equation}
Therefore,
\begin{equation}
\label{3I-G45a}
F^{ii} U_{ii}^2 \geq - K \sum_J F^{ii} U_{ii}
\geq - \frac{K}{n} \sum_{U_{ii} < 0} F^{ii} U_{ii}
\geq \frac{\varepsilon K}{n} -  C K \sum F^{ii}.
\end{equation}

Finally, note that $1 - \varphi_u d \geq \frac{1}{2}$ and $d$ small enough.
By \eqref{3I-G25b}-\eqref{3I-G45a} and
\eqref{3I-G15a} we derive
\begin{equation}
\label{3I-G45ab}
 (c_0 K^2 - C K^{\frac{3}{2}} - C K - C ) \sum F^{ii}
+ c_0 \epsilon K - C K^{\frac{1}{2}} \leq
w (\varphi_u d - 1) F^{ii} U_{ii}
\end{equation}
for some $c_0 > 0$.
This gives a bound $K \leq C$, completing our proof of the gradient estimates, provided that $f$ satisfies
\eqref{082301}.

\bigskip

\section{Global estimates for second order derivative }\label{global}
\label{S}

In this we shall derive the second order derivative estimate \eqref{hess-a10} in Theorem~\ref{mainth1}.
Motivated by \cite{LTU} and \cite{MQ} we consider the
following quantity.
Let
\begin{equation}
W (x, \xi) = \nabla_{\xi \xi} u + \chi_{\xi \xi}
   - 2 g (\xi, \nu) (\nabla_ {\xi'} \tvarphi - \nabla_{\nu_{\xi'}} u + \chi_{\xi' \nu})
\end{equation}
 for $x \in \bM$ and $\xi \in T_x M$, where $\xi' = \xi - g (\xi,\nu) \nu$,
$\nu_{\xi'} = \nabla_{\xi'} \nu$ and
$\tvarphi (x) = \varphi (x, u)$; for convenience we shall write
\[ W' (x, \xi) = 2 g (\xi, \nu) (\nabla_ {\xi'} \tvarphi
   - \nabla_{\nu_{\xi'}} u). \]

Let $\eta = \eta (u, |\nabla u|)$ be a function to be determined. We consider
\begin{equation}
 \tilde {W} =  \max_{x \in \bM}
    \max_{\xi \in T_x M, |\xi| = 1} W  e^{\eta}.
\end{equation}
Our goal is to derive a bound for $\tilde{W}$ which we shall assume to be positive.
Suppose that $\tilde{W}$ is achieved at a point
$x_0 \in \bM$ for some unit vector
$\xi \in T_{x_0} M$. We shall consider separately
two different cases: ({\bf a}) $x_0 \in M$
and ({\bf b}) $x_0 \in \partial M$.
In this section we consider case ({\bf a}) while
case ({\bf b}) will be treated in Section~\ref{B}.

Assume now that $x_0$ is an interior point. We choose smooth orthonormal
local frames $e_1, \ldots, e_n$
about $x_0 \in M$ such that $e_1= \xi$,
$\nabla_i e_j = 0$ so $\Gamma_{ij}^k = 0$
for all $1 \leq i, j, k \leq n$,
and $U_{ij} = \nabla_{ij} u + \chi_{ij}$ is diagonal  at $x_0$.
Denote $W = W (x, e_1)$ and
\[ W' = W' (x, e_1)
        = 2 g (e_1, \nu) (\nabla_ {e_1'} \varphi +
               \varphi_u \nabla_ {e'_1} u
   - \nabla_{e'_1 k} d \nabla_k u). \]
which are locally defined near $x_0$. In what follows we modify the argument in \cite{Guan14} to derive
a bound for $W$ at $x_0$; note that we can not use directly the result there though.

The function
$\log W + \eta$ attains its maximum at $x_0$
and therefore for $i = 1, \ldots, n$,
\begin{equation}
\label{hess-a30}
\frac{\nabla_i W}{W} + \nabla_i \eta = 0,
\end{equation}
and
\begin{equation}
\label{hess-a40}
\begin{aligned}
\frac{\nabla_{ii} W} {W}
   - \Big ( \frac {\nabla_i W} {W} \Big) ^2 + \nabla_{ii} \eta \leq 0.
\end{aligned}
\end{equation}


Differentiating equation \eqref{3I-10'} twice we obtain at $x_0$,
\begin{equation}
\label{hess-a65}
\begin{aligned}
F^{ii} \nabla_{11} U_{ii} + \sum F^{ij, kl} \nabla_1 U_{ij} \nabla_1 U_{kl}
 = \nabla_{11} \psi.
\end{aligned}
\end{equation}
It follows 
that
\begin{equation}
\label{hess-a68}
\begin{aligned}
F^{ii} \nabla_{ii} U_{11}
 \geq \,&  - F^{ij, kl} \nabla_1 U_{ij} \nabla_1 U_{kl}
                + \nabla_{11} \psi
                - C(1 + |\nabla u| + W)\sum F^{ii}
\end{aligned}
\end{equation}
and
\begin{equation}
\label{032601}
\begin{aligned}
F^{ii} \nabla_{ii} W'
\leq \,& 2 g (e_1, \nu) F^{ii}
      \nabla_{ii} (\varphi_u \nabla_ {e'_1} u
   - \nabla_{e'_1 k} d \nabla_k u) \\
  & + C  (1 + |\nabla u|^2 + W) \sum F^{ii} \\
\leq \, & C (1 + |\nabla u|^2 + W + Z) \sum F^{ii}
+ C,
\end{aligned}
\end{equation}
where
\[ Z = |\varphi_{uuu}| |\nabla u|^3 +
          |\varphi_{uu}| |\nabla u| W. \]
We now plug
 (\ref{hess-a68}) and \eqref{032601} in
\eqref{hess-a40} to derive
\begin{equation}
\label{032602}
W F^{ii} \nabla_{ii}  \eta \leq
E + C (1 + |\nabla u|) +  C_2 \sum F^{ii},
\end{equation}
where 
$C_2 = C (1 + |\nabla u|^2 + W + Z)$ and
\[ E = F^{ij, kl} \nabla_1 U_{ij} \nabla_1 U_{kl} +
  \frac{1}{W} \sum F^{ii} (\nabla_i W)^2. \]

As in \cite{Guan14, Guan} to estimate $E$ we follow an idea of
Urbas~\cite{Urbas02}.
Let $0 < s < 1$  to be chosen and
\[  \begin{aligned}
J \,& = \{i: U_{ii} \leq - s U_{11}\}, \\
K \,& = \{i: U_{ii} > -s  U_{11}, \; s F^{ii} \geq F^{11}\}, \\
L \,&= \{i: U_{ii} > -s  U_{11}, \; s F^{ii} < F^{11}\}.
  \end{aligned} \]
It was shown by Andrews~\cite{Andrews94}
and Gerhardt~\cite{Gerhardt96}, also earlier by  Caffarelli, Nirenberg and Spruck, that
\[ - F^{ij, kl} \nabla_1 U_{ij} \nabla_1 U_{kl}
 \geq \sum_{i \neq j} \frac{F^{ii} - F^{jj}}{U_{jj} - U_{ii}}
            (\nabla_1 U_{ij})^2. \]
By \eqref{hess-a30} and Schwarz inequality we obtain
\begin{equation}
\label{gsz-G270}
\begin{aligned}
 - F^{ij, kl} \nabla_1 U_{ij} \nabla_1 U_{kl}
 \geq \,& 2 \sum_{i \geq 2}
                \frac{F^{ii} - F^{11}}{U_{11} - U_{ii}}
                (\nabla_1 U_{i1})^2 \\
 \geq \,& \frac{2 (1-s)}{(1+s) U_{11}}
               \sum_{i \in K} F^{ii}
              (\nabla_1 U_{i1})^2 \\
 \geq \,& \frac{2 (1 - s)^2}{(1 + s) U_{11}}
             \sum_{i \in K} F^{ii}
           ((\nabla_i U_{11})^2 - C|\nabla u|^2/s)\\
\geq \,&\frac{2 (1 - s)^3}{(1 + s) U_{11}}
        \sum_{i \in K} F^{ii} (\nabla_i W)^2
             \\
  & - \frac{C}{s U_{11}} \sum_{i \in K}
       F^{ii} ( |\nabla u|^2 + (\nabla_i W')^2) .
\end{aligned}
\end{equation}
By straightforward calculations,
\begin{equation}
\label{gsz-G275}
|\nabla W'|^2 \leq
  C (W^2 + |\nabla u|^2 + \varphi_{uu} |\nabla u|^4).
\end{equation}

Next, we may assume
$|W'| \leq s U_{11}$ at $x_0$ and
fix $s = 1/9$ so that
\begin{equation}
\label{032606}
(1 - s) U_{11} \leq W \leq (1 + s) U_{11}
\end{equation}
and
\begin{equation}
\label{032605}
\frac{2 (1 - s)^3 W}{(1 + s) U_{11}} \geq
\frac {2 (1 - s)^4}{(1 + s)} \geq 1.
 \end{equation}
By \eqref{gsz-G270}-\eqref{032605} we obtain
\begin{equation}
\label{032604}
\begin{aligned}
E \leq\,& \frac{1}{W}
        \sum_{i \in J \cup L} F^{ii} (\nabla_i W)^2
+ C_3 \sum_{i \in K} F^{ii},
           \end{aligned}
\end{equation}
where
\[ C_3 = C  W + \frac{|\nabla u|^2
           (1 + |\varphi_{uu}| |\nabla u|^2)}{W}. \]

Following \cite{Guan} we take
$\phi(t) = -\log (1 - \gamma b t)$ and
\[ \eta = \phi (|\nabla u|^2) + a (v - u), \]
where $a > 0$ and $\gamma \in (0, 1/2]$ are constant
to be determined, $b = 1/2 b_1$ and
\[ b_1 = \max_{\bM} (1 + |\nabla u|^2). \]
So
\[ \frac{\gamma b}{2} \leq \phi' (|\nabla u|^2)
= \sqrt{ \phi '' (|\nabla u|^2) }
 = \frac{\gamma b} {1 - bt} \leq 2 \gamma b. \]
We calculate
\begin{equation}
 \begin{aligned}
 \nabla_i \eta
   = \,& 2 \phi' \nabla_{k} u \nabla_{ik} u + a \nabla_i (v - u)
   = 2 \phi' (U_{ii} \nabla_i u )
        + a \nabla_i (v - u),
 \end{aligned}
  \end{equation}
\begin{equation}
\begin{aligned}
  \nabla_{ii} \eta
   = \,&  2 \phi' (\nabla_{ik} u \nabla_{ik} u
          + \nabla_{k} u \nabla_{iik} u)
          + 4 \phi'' (\nabla_{k} u \nabla_{ik} u)^2 
          + a \nabla_{ii} (v - u).
\end{aligned}
 \end{equation}
Therefore,
\begin{equation}
\label{hess-a272}
\begin{aligned}
  \sum_{i \in J \cup L} F^{ii} (\nabla_i \eta)^2
  \leq \, &
   C a^2 |\nabla u|^2 \sum_{i \in J \cup L} F^{ii}
  + C (\phi')^2 |\nabla u|^2 F^{ii} U_{ii}^2
\end{aligned}
\end{equation}
and
 \begin{equation}
\label{hess-a273}
 \begin{aligned}
 F^{ii} \nabla_{ii} \eta
\geq \,&  \phi'  F^{ii} U_{ii}^2 + 4 \phi'' F^{ii} (\nabla_{k} u \nabla_{ik} u)^2
       + a F^{ii} \nabla_{ii} (v - u) \\
 & - C\phi'|\nabla u|^2  \sum F^{ii}-C\phi'|\nabla u|.
\end{aligned}
\end{equation}
By \eqref{hess-a30}, \eqref{032602}, \eqref{032604},
\eqref{hess-a272} and \eqref{hess-a273}, we derive
\begin{equation}
\label{032608}
\begin{aligned}
\phi'  F^{ii} U_{ii}^2 + a F^{ii} \nabla_{ii} (v - u)
 \leq \,&
C a^2 |\nabla u|^2 \sum_{i \in J \cup L} F^{ii}
+ C (\phi')^2 |\nabla u|^2 F^{ii} U_{ii}^2 \\
& + C (1 + |\varphi_{uu}| |\nabla u|
    +  |\varphi_{uuu}| |\nabla u|^3 W^{-1})  \sum F^{ii} \\
 &
+ \frac{C |\nabla u|^2 (1 + |\varphi_{uu}| |\nabla u|^2)}{W^2} \sum F^{ii} - \frac{\nabla_{11} \psi}{W} + C.
\end{aligned}
\end{equation}

Next, note that
\begin{equation}
\label{hess-a274}
 F^{ii} U_{ii}^2 \geq  F^{11} U_{11}^2
      + \sum_{i \in J} F^{ii} U_{ii}^2
   \geq F^{11} U_{11}^2
       + s^2 U_{11}^2 \sum_{i \in J} F^{ii}.
 \end{equation}
We now fix $\gamma \leq 1/16 C$ so that
\[ \phi' - C (\phi')^2 |\nabla u|^2 \geq
\frac{\gamma b}{2} - 4 C \gamma^2 b
=  \frac{\gamma b (1- 8 C \gamma)}{2}
\geq \frac{\gamma b}{4}. \]
By Theorem~2.17 in \cite{Guan14} we may fix
$a \geq A (1 + |\varphi_{uu}| |\nabla u|)$ for $A$ sufficiently large to obtain
\[ a W \leq C (|\varphi_{uuu}| |\nabla u|^3
      + |\nabla^2 \psi|) \]
or
\[ U_{11} (x_0) \leq C a (1 + |\nabla u|^2)
\leq C (1 + |\nabla u|^2
    + |\varphi_{uu}| |\nabla u|^3). \]
Consequently,
\[ W \leq C (1 + |\nabla u|^2)
  + C (|\varphi_{uu}| + |\varphi_{uuu}|) |\nabla u|^3. \]
 where $ C $ is independent of
$|\nabla u|_ {C^0 (\bM)}$.
The proof of Theorem~\ref{mainth1} is complete.

\bigskip

\section{Second derivative estimates on boundary}
\label{B}

\medskip

In this section we consider case ({\bf b})
$x_0 \in \partial M$. First assume that
$0 < |\xi'| < 1$ and let
$\tau = |\xi'|^{-1} \xi'$. We have
\[ \begin{aligned}
\nabla_{\xi \xi} u = \,& \nabla_{\xi' \xi'} u
   + 2 g(\xi, \nu)  \nabla_{\xi' \nu} u
   +g(\xi, \nu)^2 \nabla_{\nu \nu} u \\
    = \,& |\xi'|^2 \nabla_{\tau \tau} u
   + 2 g(\xi, \nu)  (\nabla_{\xi'} (\nabla_{\nu} u)
-  \nabla_{\nu_{\xi'}} u)
   + (1 - |\xi'|^2) \nabla_{\nu \nu} u \\
   = \,& |\xi'|^2 \nabla_{\tau \tau} u
   + 2 g(\xi, \nu)  (\nabla_{\xi'} \tvarphi
-  \nabla_{\nu_{\xi'}} u)
   + (1 - |\xi'|^2) \nabla_{\nu \nu} u.
\end{aligned}\]
It follows that
 \[ \begin{aligned}
W (x_0, \xi) = \,& |\xi'|^2 W (x_0, \tau)
   + (1- |\xi'|^2) W (x_0, \nu) \\
 \leq \,& |\xi'|^2 W (x_0, \xi)
   + (1- |\xi'|^2) W (x_0, \nu)
 \end{aligned} \]
which implies
\begin{equation}
W (x_0, \xi) \leq  W (x_0, \nu).
\end{equation}
Consequently, we only need to considier the following
two cases: ({\bf i}) $|\xi'| = 1$
and ({\bf ii}) $|\xi'| = 0$.

Case ({\bf i}) $|\xi'| = 1$.  So $\xi$ is tangential to $\partial M$ at $x_0$. We choose smooth
orthonormal local frames $e_1, \ldots, e_n$ around $x_0$
such that $e_n = \nu$ along $\partial M$
and $e_1 = \xi$ at $x_0$.
Write $U_{11} = \nabla_{11} u + \chi_{11}$.
At $x_0$ we have
\begin{equation}
\label{9-26-7}
\begin{aligned}
0 \geq \,&  \nabla_n W + W \nabla_n \eta
   \geq  \nabla_n U_{11} - \nabla _n W'
           + W \nabla_n \eta.
\end{aligned}
\end{equation}

Since $g (\xi, \nu) = 0$, we have $W' = 0$ and
\begin{equation}
\label{9-26-85}
    \nabla_n W' \geq - C (1 + |\nabla u|) .
\end{equation}
By the boundary condition~\eqref{021906}
 we have for $k < n$,
\begin{equation}
\label{9-26-35}
 \nabla_{kn} u = \nabla_k (\nabla_n u)
  - \nabla_{\nabla_k e_n} u
= \nabla_k \tvarphi + b_{kl} \nabla_l u,
\end{equation}
where $\{b_{kl}\}$ denotes the second fundamental form of $\partial M$, and
\begin{equation}
\label{9-26-95}
\begin{aligned}
\nabla_{11n} u
      \,& = \nabla_{11} (\nabla_n u)
               - 2 g(\nabla_{e_1} e_n, e_1) \nabla_{11} u
               - \nabla_{\nabla_{11} e_n} u \\
\,& \geq \nabla_{11} \tvarphi + 2 b_{11} \nabla_{11} u
              - C |\nabla u| \\
\,& \geq (\varphi_u + 2 b_{11}) \nabla_{11} u
              - C (1 + |\nabla u|^2).
\end{aligned}
\end{equation}

Next, note that $\nabla_{1k} u (x_0) = 0$ for all
$k > 1$. Indeed, 
let
\[ e_{\theta} = e_1 \cos \theta + e_k \sin \theta. \]
We have
\[ \nabla_{e_{\theta} e_{\theta}} u =
\nabla_{11} u \cos^2 \theta
 + 2 \nabla_{1k} u \cos \theta \sin \theta
 + \nabla_{kk} u \sin^2 \theta. \]
It follows that
\[ \frac{\partial}{\partial \theta}
          \nabla_{e_{\theta} e_{\theta}} u
   = (\nabla_{kk} u - \nabla_{11} u)
          \sin 2 \theta +  \nabla_{1k} u \cos 2 \theta. \]
and
\[2 \nabla_{1k} u = \frac{\partial}{\partial \theta} \Big|_{\theta = 0} \nabla_{e_{\theta} e_{\theta}} u = 0. \]
On the other hand, $g (\nabla_n e_1, e_1) = 0$ since
$|e_1|= 1$.
By \eqref{9-26-95} it follows that
\begin{equation}
\label{9-26-45}
\begin{aligned}
\nabla_n U_{11}
 = \,& \nabla_n (\nabla_{11} u)
            + \nabla_n \chi_{11} \\
 = \,&  \nabla_{n11} u  + \nabla_n \chi_{11} \\
 = \,& \nabla_{11n} u + R_{l11n} \nabla_l u
            + \nabla_n \chi_{11} \\
\geq \,&(\varphi_u + 2 b_{11}) \nabla_{11} u
              - C (1 + |\nabla u|^2) \\
\geq \,&(\varphi_u + 2 b_{11}) W
              - C (1 + |\nabla u|^2).
\end{aligned}
\end{equation}

From Section~\ref{S},
\[ \eta = \phi (|\nabla u|^2) + a (v - u), \]
where $a$ is a sufficiently large positive constant,
$v$ as in \eqref{032201},
and $\phi$ is a nondecreasing function satisfying
$t \phi' (t) \leq 1$.
By \eqref{9-26-35},
\begin{equation}
\label{9-26-75}
\begin{aligned}
    \nabla_n \eta
     = \,& 2 \phi' \nabla_k u \nabla_{nk} u
             + a \nabla_n (v - u) \\
\geq \,& \epsilon_0 a + 2 \phi' \tvarphi \nabla_{nn} u
           - C \phi' (1 + |\nabla u|^2).
\end{aligned}
\end{equation}

Finally, it follows from \eqref{9-26-7}, \eqref{9-26-85},
\eqref{9-26-45} and \eqref{9-26-75} that
\begin{equation}
\label{9-26-125}
\begin{aligned}
(\epsilon_0 a + 2 \phi' \varphi \nabla_{nn} u
   + \varphi_u + 2 b_{11} - C \phi' |\nabla u|^2) W
\leq C (1 + |\nabla u|^2).
\end{aligned}
\end{equation}
By
\eqref{9-26-115} below we derive
\begin{equation}
\label{9-26-105}
W \leq C (1 + |\nabla u|)
\end{equation}
provided that $a$ is sufficiently large, independent of $|u|_{C^1 (\bM)}$, so that
\[ \epsilon_0 a + 2 \phi' \varphi \nabla_{nn} u
   + \varphi_u + 2 b_{11} - C \phi' |\nabla u|^2
 \geq 1. \]

In order to derive a bound for $W$ in both cases ({\bf i}) and ({\bf ii}), it is therefore enough to establish the double normal derivative estimate
\begin{equation}
\label{9-26-115}
|\nabla_{\nu \nu} u|  \leq C (1 + |\nabla u|)
\;\; \mbox{on $\partial M$.}
\end{equation}
The rest of this section is devoted to this estimate.

Consider an arbitrary point $x_0 \in\partial M$. We choose
smooth orthonormal local frames
$e_1, \ldots, e_n$ around $x_0$ as before such that $e_n = \nu$ along $\partial M$.
We shall use $\rho(x)$ to denote the distance function from $x$ to $x_0$,
\begin{equation}
\rho(x) \equiv \mbox{dist}_{M} (x, x_0),
\end{equation}
and let
$M_{\delta} = \{x \in M: \rho (x) < \delta\}$.
Note that $\nabla_{ij} \rho^2 (x_0)=2\delta_{ij}$. We may assume
\begin{equation}
 \{\delta_{ij}\} \leq \{\nabla_ {ij} \rho^2\}
                       \leq 3 \{\delta_{ij}\}
\;\; \mbox{in $M_{\delta}$}.
\end{equation}

Differentiating equation~\eqref{3I-10'} we obtain
near $x_0$,
\begin{equation}
\label{hess-a60'}
\begin{aligned}
F^{ij} \nabla_{ij} (\nabla_k u)
= \, & F^{ij} (\nabla_{ijk}  u
      + \Gamma_{jk}^l \nabla_{il}  u
      + \Gamma_{ik}^l \nabla_{jl}  u) \\
= \,& F^{ij} \nabla_{kij}  u + F^{ij} R_{ijkl} \nabla_l u
        + 2 F^{ij} \Gamma_{jk}^l \nabla_{il}  u \\
= \, &\nabla_k \psi + F^{ij} R_{ijkl} \nabla_l u
        + 2 F^{ij} \Gamma_{jk}^l \nabla_{il}  u,
\end{aligned}
\end{equation}
where $R_{ijkl}$ and $\Gamma_{jk}^l$
denote the Riemannian curvature tensor and Christoffel
symbols, respectively.

We now make use of the assumptions that $M$ is locally conformally flat 
near $\partial M$ and $\partial M$ is umbilic. 
It follows that when $\delta$ is sufficiently
small,
\begin{equation}
\Gamma_{jn}^l = g (e_l, \nabla_j e_n) = - \kappa_d \delta_{jl}
   \;\; \mbox{in $M_{\delta}$},
\end{equation}
where $\kappa_d (x)$ denotes the principal curvature of the level hypersurface of $d$ passing through
$x \in M_{\delta}$.
Therefore,
\begin{equation}
\label{hess-a610}
\begin{aligned}
F^{ij} \nabla_{ij} (\nabla_n u - \tvarphi)
= \,& 2 F^{ij} \Gamma_{jn}^l \nabla_{il} u
- \varphi_u F^{ij} \nabla_{ij}  u + Q \\
=  \,&  - (2 \kappa_d + \varphi_u) F^{ij} U_{ij} + Q, 
\end{aligned}
 \end{equation}
where
\[ |Q| \leq C (1 + |\nabla u| ) \sum F^{ii} + C, \]
 and
\begin{equation}
\label{hess-a620}
\begin{aligned}
F^{ij} \nabla_{ij} (\nabla_n u - \tvarphi)^2
= \,& 2 (\nabla_n u - \tvarphi)
          F^{ij} \nabla_{ij} (\nabla_n u - \tvarphi) \\
      & + 2 F^{ij} \nabla_i (\nabla_n u - \tvarphi)
                        \nabla_j (\nabla_n u - \tvarphi) \\
\geq \,& F^{ij} U_{in} U_{nj}  - C |F^{ij} U_{ij}|
              - C (1 + |\nabla u| ) \sum F^{ii}.
\end{aligned}
 \end{equation}

We now construct a barrier function as follows. Let
\begin{equation}
H = 
    A_1 \underline{u}  - A_2 (d - N d^2) - A_3 \rho^2,
\end{equation}
where $A_1$, $A_2$, $A_3$ and $N$ are positive constants to be chosen later; we shall
assume $2 N \delta \leq 1$ by fixing $\delta$ small after $N$ is determined. We first calculate
\begin{equation}
\label{032306}
\begin{aligned}
F^{ij} \nabla_{ij} (d - N d^2)
    = \, & (1 - 2 N d) F^{ij} \nabla_{ij} d
              - 2 N F^{ij} \nabla_i d \nabla_j d \\
\leq \,& - 2 N F^{ij} \nabla_i d \nabla_j d
            + C \sum F^{ii}.
\end{aligned}
\end{equation}

At a fixed point $x \in \ol M_{\delta}$,
let $\mu = \lambda [\nabla^2 \ul u]$ and
$\lambda = \lambda [U]$ denote the eigenvalues of
$\nabla^2 \ul u$ and $U$, respectively.
By assumption~\eqref{032301} and
Theorem 2.17 in \cite{Guan14}, there are positive constants $R$ and $\varepsilon$ depending
 such that
\begin{equation}
\label{GX-B150}
 \sum f_i (\lambda) (\mu_i - \lambda_i) \geq
    \varepsilon \Big(1 + \sum f_i\Big)
\end{equation}
provided that $|\lambda| \geq R$.

We first assume $|\lambda| \leq R$.
In this case, there are uniform bounds
\begin{equation}
0 < \alpha \leq f_i (\lambda) \leq \beta,
\end{equation}
with $\alpha$, $\beta$ depending on $R$. Therefore,
$F^{ij} \nabla_i d \nabla_j d \geq \alpha$ since
$|\nabla d| \equiv 1$. We now fix $N$ sufficiently large
and then  $\delta > 0$ such that
$2 N \delta \leq 1$ and, by \eqref{032306},
\begin{equation}
\label{032306a}
\begin{aligned}
F^{ij} \nabla_{ij} (d - N d^2)
\leq \,& - \alpha N \sum F^{ii} \leq - n \alpha^2 N.
\end{aligned}
\end{equation}

It is clear that $H(0) = 0$ and
\[ H \leq - A_3 \rho^2 \;\; \mbox{in $M_{\delta}$}. \]
Since $\nabla_n u - \tvarphi$ = 0 on $\partial M$,
we may fix $A_3$ and then $A_2$ large such that
\begin{equation}
 |\nabla_n u - \tvarphi| \leq A_3 \rho^2 \leq - H \;\;
\mbox{on $\partial M_{\delta}$},
\end{equation}
and
\begin{equation}
\label{091301}
\begin{aligned}
 F^{ij} \nabla_{ij} (H \pm (\nabla_{\nu} u - \varphi))
\geq \,& 
   (A_2 \alpha N - 3 A_3 - C - C R) \sum F^{ii}
 - C \geq 0.
\end{aligned}
\end{equation}
Here we have used the fact that
\begin{equation}
\label{GX-B160a}
  F^{ij} \nabla_{ij} \ul u
\geq F (\nabla^2 \ul u) - F (U) + \sum F^{ij} U_{ij}
\geq 0 
\end{equation}
by the concavity of $f$ and assumptions
\eqref{gx-I110a} and \eqref{032301}.

Suppose now that $|\lambda| > R$. We note that
by assumption~\eqref{gx-I110a} and the concavity of $f$,
\[ 0 \leq F^{ij} U_{ij} \leq F (U) - F (g) + \sum F^{ii}
   \leq \sum F^{ii} + C. \]
By \eqref{GX-B150} we may finally fix $A_1$ such that
\begin{equation}
\label{091301a}
\begin{aligned}
 F^{ij} \nabla_{ij} (H \pm (\nabla_{\nu} u - \varphi))
\geq \,& A_1 F^{ij} \nabla_{ij} \ul u
  - C (A_2 + 3 A_3 + 1) \sum F^{ii} - C \geq 0.
\end{aligned}
\end{equation}

Consequently, by the maximum principle we obtain
\[ H \pm (\nabla_{\nu} u - \varphi )\leq 0 \;\;
\mbox{in $M_{\delta}$} \]
and therefore,
\begin{equation}
\label{GX-B160}
|\nabla_{nn} u| \leq \nabla_n H + |\nabla_n \tvarphi|
\leq C (1 + |\nabla u|).
\end{equation}
This completes the proof of Theorem~\ref{3I-thm2}.


\end{document}